\documentclass[a4paper,10pt,final]{amsart}
\usepackage{amssymb}
\usepackage{amsmath}

\renewcommand{\phi}[0]{\varphi}
\renewcommand{\theta}[0]{\vartheta}

\newcommand{\F}{\mathbb{F}}

\DeclareMathOperator{\Irr}{Irr}

\newtheorem{dummy}{Dummy}

\numberwithin{dummy}{section}
\numberwithin{equation}{section}

\newtheorem{lemma}[dummy]{Lemma}
\newtheorem{theorem}[dummy]{Theorem}
\newtheorem{prop}[dummy]{Proposition}
\newtheorem{cor}[dummy]{Corollary}

\theoremstyle{definition}

\theoremstyle{remark}
\newtheorem{rem}[dummy]{Remark}

\begin{document}
\bibliographystyle{amsalpha}

\author{S.~Mattarei}
\email{mattarei@science.unitn.it}
\urladdr{http://www-math.science.unitn.it/\~{ }mattarei/}
\address{Dipartimento di Matematica\\
  Universit\`a degli Studi di Trento\\
  via Sommarive 14\\
  I-38050 Povo (Trento)\\
  Italy}
\title[Nonsingular derivations]{The orders of nonsingular derivations of Lie algebras of characteristic two}
\thanks{The  author  is grateful  to  Ministero dell'Istruzione e dell'Universit\`a,
  Italy,  for  financial  support to  the
  project ``Graded Lie algebras  and pro-$p$-groups of finite width''.}

\begin{abstract}
Nonsingular derivations of modular Lie algebras which have finite multiplicative order play a role
in the coclass theory for pro-$p$ groups and Lie algebras.
A study of the set $\mathcal{N}_p$ of positive integers which occur as
orders of nonsingular derivations of finite-dimensional non-nilpotent Lie algebras
of characteristic $p>0$ was initiated by Shalev and continued by the present author.
In this paper we continue this study in the case of characteristic two.
Among other results, we prove that any divisor $n$ of $2^k-1$ with $n^4>(2^k-n)^{3}$
belongs to $\mathcal{N}_2$.
Our methods consist of elementary arguments with polynomials over finite fields
and a little character theory of finite groups.
\end{abstract}



\maketitle

\thispagestyle{empty}

\section{Introduction}\label{sec:intro}

A classical result of Jacobson (proved in~\cite{Jacobson:autom-and-deriv} and also listed as Problem~9 of Chapter~II
in~\cite[p.~54]{Jac:Lie_algebras})
asserts that if a finite-dimensional Lie algebra $L$
over a field of characteristic zero admits a derivation $D$ which is nonsingular
(that is, injective), then $L$ is nilpotent.
Jacobson also proved that the analogous result holds in prime characteristic
provided $L$ is restricted.
In absence of the restrictedness hypothesis some other assumption on $L$ or $D$ is necessary
to preserve the conclusion that $L$ is nilpotent.

A striking instance of such a result occurs in the effective proof given by Shalev in~\cite{Sha:coclass}
of the strongest of the coclass conjectures of Leedham-Green and Newman
for pro-$p$ groups~\cite{L-GN}.
A simple but crucial step in Shalev's proof is the fact that
a finite-dimensional
Lie algebra over a field of characteristic $p>0$
having a derivation $D$ with $D^{p-1}=1$,
that is, nonsingular and with all eigenvalues in the prime field,
must be nilpotent~\cite[Proposition 5.2.8]{L-GMcKay}.

Because of the importance of the above fact it would be interesting to know to what extent the hypothesis
on the order of $D$ can be weakened.
There are at least two natural ways in which this hypothesis may be weakened,
the first one being imposing that $D^{p^k-1}=1$
for some positive integer $k$.
This hypothesis is insufficient as soon as $k>1$,
because for all $k>1$ there exist even simple finite-dimensional Lie algebras
of characteristic $p$ which have
nonsingular derivations of order $p^k-1$, namely, certain algebras
discovered by Albert and Frank~\cite{AF} in the fifties.
These Lie algebras, which belong to the larger class of {\em Block algebras}
and are usually denoted by $H(2:\mathbf{n};\Phi(\tau))^{(2)}$ as in~\cite{Strade:book},
were employed in \cite{Sha:max} to construct the first examples
of non-soluble modular graded Lie algebras of maximal class,
thus disproving the analogues of Conjectures~C and~D of \cite{L-GN} for modular graded Lie algebras.
They later turned out to be the building blocks for the construction of all graded Lie algebras of maximal class
(generated by their homogeneous component of degree one)
over fields of odd characteristic~\cite{CMN,CN}.
We should also mention in this connection that Benkart, Kostrikin and Kuznetsov have determined in~\cite{BKK}
all finite-dimensional Lie algebras over an algebraically closed field of characteristic $p>7$
which admit a nonsingular derivation.
We refer to the Introduction of~\cite{CaMa:Hamiltonian} for a broader discussion and references on these
(and related) topics.

A second way of weakening the assumption that $D^{p-1}=1$ is imposing an upper bound on the order of
the derivation $D$.
In this direction, Shalev proved in~\cite{Sha:nonsing-der} that a finite-dimensional
modular Lie algebra with a nonsingular derivation
of order less than $p^2-1$ must be nilpotent, where $p$ is the characteristic of the underlying field,
which we assume without mention from now on.
This result is best possible because of the Lie algebras of Albert and Frank mentioned in the previous paragraph.

Shalev suggested in~\cite{Sha:nonsing-der} the more general problem of studying
the set $\mathcal{N}_p$ of positive integers which occur as the orders of nonsingular
derivations of finite-dimensional non-nilpotent Lie algebras of prime characteristic $p$.
Thus, we know that $\mathcal{N}_p$ contains all numbers of the form $p^k-1$, for all $k\ge 2$,
and Shalev proved that $p^2-1$ is the smallest element in $\mathcal{N}_p$.
We extended Shalev's result in~\cite{Mat:nonsing-der} by proving that for $p>3$ the only
numbers in $\mathcal{N}_p$ which are less than $p^3-1$ are multiples of $p^2-1$.

The key for proving these and similar results is the following translation of the problem
into one formulated entirely in terms of finite fields.
A positive integer $n$ belongs to $\mathcal{N}_p$ if and only if there is an element $\alpha\in \bar\F_p$
(the algebraic closure of the field of $p$ elements $\F_p$) such that
$(\alpha+\lambda)^n=1$ for all $\lambda\in\F_p$.
The necessity of the condition was proved in~\cite{Sha:nonsing-der} by means of the Engel-Jacobson theorem,
and the sufficiency in~\cite{Mat:nonsing-der} by means of an explicit construction.
We comment more on this characterization of $\mathcal{N}_p$ and a further one at the beginning of Section~\ref{sec:even}.

In the present paper we continue the study of $\mathcal{N}_p$ which was initiated
in~\cite{Sha:nonsing-der} and~\cite{Mat:nonsing-der},
with special emphasis on the case where $p=2$, because of reasons which we explain below.
The characterization of $\mathcal{N}_p$ recalled in the previous paragraph allows us to
forget about the origin of the problem in the theory of Lie algebras.
Thus, no knowledge of Lie algebras is necessary to understand the paper beyond this Introduction.
We briefly describe the contents of the paper. 

It is easy to see that $\mathcal{N}_p$ contains all multiples of its elements.
Shalev exhibited in~\cite{Sha:nonsing-der} some elements of $\mathcal{N}_p$ which are not multiples
of $p^k-1$ for any $k>1$, namely, $(p^p-1)/(p-1)$ for $p$ odd, and $73$ for $p=2$.
In Section~\ref{sec:even} we generalize this latter example to show that
$(q^3-1)/(q-1)=q^2+q+1\in\mathcal{N}_2$, for all $q=2^s$ with $s\ge 1$.
The argument lends itself to a further generalization which might be of independent interest,
on the period of a polynomial of the form $cx^{q+1}+dx^q-ax-b\in\F_q[x]$, where $q=p^s$ and $p$ is any prime.

In Section~\ref{sec:Frobenius} we
show the abundance of elements of $\mathcal{N}_2$
by proving that all divisors of $2^k-1$ which are ``large enough'' belong to $\mathcal{N}_2$.
More precisely, Theorem~\ref{thm:U-big}, the main result of this paper, states that
if $n$ is a divisor of $2^k-1$ such that $n^4>(2^k-n)^3$, then $n\in\mathcal{N}_2$.
It follows, in particular, that
$(q^t-1)/(q-1)\in\mathcal{N}_2$
if $q=2^s$ with $s\ge 1$ and $t>3$.
(The case $t=3$ escapes Theorem~\ref{thm:U-big} but is dealt with directly in Proposition~\ref{prop:2^2s+};
the case $t=4$ admits also an independent proof, as in Proposition~\ref{prop:2^3s+}.)
The proof of Theorem~\ref{thm:U-big} is based on the character theory of a certain Frobenius group,
but the method does not appear to extend to characteristics higher than two.
These require more sophisticated tools, and we plan to deal with them in a future paper.
We conclude Section~\ref{sec:Frobenius} with a discussion of the set of multiples of numbers of this form.

In Section~\ref{sec:counting} we introduce a problem which extends the mere determination of $\mathcal{N}_p$.
For any positive integer $n$ prime to $p$ let $E_p(n)$ denote the set of elements
$\alpha\in \bar\F_p$
such that
$(\alpha+\lambda)^n=1$ for all $\lambda\in\F_p$.
Then $n\in\mathcal{N}_p$ occurs exactly when $E_p(n)$ is not empty.
We hope that a study of $E_p(n)$, or at least of its cardinality,
may shed some light on the structure of the set $\mathcal{N}_p$,
even though it is not presently clear what significance this additional information may have
for nonsingular derivations of Lie algebras.
Besides its intrinsic general interest, this study has also computational motivations which we explain in Remarks~\ref{rem:remainder}
and~\ref{rem:constructive}.
After the technical Lemma~\ref{lemma:remainder}, which gives an alternative and somehow more convenient
way of determining $E_p(n)$ than that suggested by its definition,
in Proposition~\ref{prop:e-values} we explicitly determine the set $E_2(n)$ for $n=(2^{3s}-1)/(2^s-1)$
as the set of roots of a certain polynomial.
A similar approach to $n=(2^{4s}-1)/(2^s-1)$ only produces a direct proof that $E_2(n)$ is not empty,
and thus $n\in\mathcal{N}_2$, in Proposition~\ref{prop:2^3s+}.
Although this statement is also a consequence of Theorem~\ref{thm:U-big}, as noted in the previous paragraph,
the proof of Proposition~\ref{prop:2^3s+} has a constructive content explained in Remark~\ref{rem:constructive}.
We conclude this section by producing a complete description of $E_p(n)$ for $n=(p^p-1)/(p-1)$,
which is the example of nontrivial element of $\mathcal{N}_p$ for $p$ odd exhibited by Shalev in~\cite{Sha:nonsing-der}.

\section{Some more numbers in $\mathcal{N}_2$}\label{sec:even}

We quote one of the main results of~\cite{Mat:nonsing-der},
which appears there as Corollary~2.3.

\begin{theorem}\label{thm:n-equiv}
Let $p$ be a prime number and let $n$ be a positive integer.
The following conditions on $n$ are equivalent:
\begin{enumerate}
\item
there exists a finite-dimensional non-nilpotent Lie algebra of characteristic $p$
with a nonsingular derivation of order $n$;
\item
there exists an element $\alpha\in\bar \F_p$ such that
$(\alpha+\lambda)^n=1$ for all $\lambda\in \F_p$;
\item
there exists an element $c\in\bar \F_p^{\ast}$ such that
$x^p-x-c$ divides $x^n-1$ as elements of the polynomial ring
$\bar \F_p[x]$.
\end{enumerate}
\end{theorem}

The version given in~\cite{Mat:nonsing-der} has the additional hypothesis that $n$ is prime to $p$.
This is superfluous since each of the three conditions holds for $n$ if and only if it does
for its $p'$-part.
(For the first condition this is shown in~\cite{Mat:nonsing-der} or~\cite{Sha:nonsing-der}.)
As we have recalled in the Introduction, the contribution of~\cite{Mat:nonsing-der} to this result
is a proof that the second condition implies the first one, while the converse
had already been proved in~\cite{Sha:nonsing-der}.
Although Lie algebras provide the motivation for this study, we will not use them in this paper,
but only investigate the set $\mathcal{N}_p$ of numbers which satisfy the second condition.

Before proceeding with the study of $\mathcal{N}_2$ in this and the next section,
we briefly elaborate on the second and third condition in Theorem~\ref{thm:n-equiv}.
For $n$ prime to $p$ let $U_n$
denote the (unique) subgroup of $\bar\F_p^\ast$ of order $n$.
Then condition~(2) has the geometric interpretation that $U_n$ contains an affine $\F_p$-line
(that is, a one-dimensional affine $\F_p$-subspace of $\bar\F_p$) with direction $1$.
(In Section~\ref{sec:counting} we address the more general problem of
determining the number of such lines.)
It follows that $\mathcal{N}_p$ contains all numbers of the form $p^r-1$ with $r>1$.
In fact, as we have mentioned in the Introduction,
for each of these numbers there exists even a simple Lie algebra of characteristic $p$
with a nonsingular derivation of order $n$.
Since $\mathcal{N}_p$ is evidently closed under taking multiples,
it contains all multiples of numbers of the form $p^r-1$ with $r>1$.
We may call these the {\em trivial} elements of $\mathcal{N}_p$.
Condition~(2) is clearly equivalent to the $p$ polynomials $(x+\lambda)^n-1$ for $\lambda\in\F_p$
having a nontrivial common factor, and this suggests an algorithm to check whether a specific
number $n$ belongs to $\mathcal{N}_p$.
In Lemma~\ref{lemma:remainder} we present a different algorithm, based on condition~(3),
which has some advantages over the former, as we discuss in Remark~\ref{rem:remainder}.
It appears that condition~(3) is also more suited than condition~(2)
to dealing theoretically with small values of $n$.
In fact, the complete description of all elements of $\mathcal{N}_p$ smaller than $p^3$,
which we obtained in~\cite{Mat:nonsing-der} extending a result of~\cite{Sha:nonsing-der}
(see the Introduction), is based on condition~(3).

In~\cite{Sha:nonsing-der} Shalev exhibited one non-trivial element of $\mathcal{N}_p$ for each prime $p$.
Since the period of the polynomial $x^p-x-1$ divides $(p^p-1)/(p-1)$
(see~\cite[Example~2.6]{Sha:nonsing-der}),
this number belongs to $\mathcal{N}_p$,
according to condition~(3).
It is easy to see that $(p^p-1)/(p-1)$ is prime to any number of the form
$p^r-1$ with $p\nmid r$.
In particular, for $p$ odd it is not divisible by any number of the form $p^r-1$ with $r>1$,
and hence it is a non-trivial element of $\mathcal{N}_p$.
For the remaining case $p=2$ Shalev showed in~\cite[Example~2.5]{Sha:nonsing-der} that $73\in\mathcal{N}_2$,
quoting from the table in~\cite[p.~378]{LN} the fact that the (irreducible) polynomial
$x^9+x+1\in\F_2[x]$ has period $73=(8^3-1)/(8-1)$.
Thus, if $\alpha$ is any root of the polynomial (in $\F_{512}$) then $\alpha^{73}=1$, and also
$(\alpha+1)^{73}=1$ since $\alpha+1=\alpha^9$,
proving that $73\in\mathcal{N}_2$ according to condition~(2).

The particular form of the polynomial used by Shalev suggests a direct computation of its period,
and the following extension of Shalev's example.
In the paper we will often use the standard notation $q$ for $2^s$ or $2^k$, depending on the context.
In doing this, we will implicitly assume that $q$ is the cardinality of a field, and thus that $s,k\ge 1$.

\begin{prop}\label{prop:2^2s+}
If $q=2^s$, then the number
$(q^3-1)/(q-1)=q^2+q+1$
belongs to $\mathcal{N}_2$.
\end{prop}

\begin{proof}
Consider the polynomial $x^{q+1}+x+1$, where $q=2^s$.
Its period  divides $(q^3-1)/(q-1)=q^2+q+1$, because
\[
x^{q^2+q+1}=(x^{q+1})^qx\equiv (x+1)^qx=x^{q+1}+x\equiv 1\pmod{x^{q+1}+x+1}.
\]
In fact, its period equals $q^2+q+1$ according to~\cite[Proposition~2.3]{Bartholdy}, but we do not need that fact here.
Therefore, each root $\alpha$ of $x^{q+1}+x+1$
has (multiplicative) order dividing $q^2+q+1$.
Furthermore, since $\alpha+1=\alpha^{q+1}$ is a power of $\alpha$, its order
also divides $q^2+q+1$.
(In fact, $\alpha+1$ has the same order as $\alpha$,
because $(q+1,q^2+q+1)=1$, but this is not essential in the argument.)
We conclude that $q^2+q+1\in\mathcal{N}_2$.

Although the proof is complete, in view of a generalization it is instructive to
do the crucial computation again, in a more conceptual way.
Congruences between polynomials will now denote equality of their images in the quotient ring
$\F_2[x]/(x^{q+1}+x+1)$.
In particular, we have
$x^q\equiv 1+1/x=(x+1)/x$, because the image of $x$ is invertible in
the quotient ring.
Since taking $q$th powers is a ring automorphism of
$\F_2[x]/(x^{q+1}+x+1)$,
we have
\[
x^{q^2}\equiv \left(\frac{x+1}{x}\right)^q=
\frac{x^q+1}{x^q}\equiv
\frac{1/x}{1+1/x}=\frac{1}{x+1},
\]
and we conclude that
$x^{1+q+q^2}=x\cdot x^q\cdot x^{q^2}\equiv x\cdot\frac{x+1}{x}\cdot\frac{1}{x+1}\equiv 1$.
\end{proof}

Note that the polynomial $x^{q+1}+x+1$ is almost never irreducible.
In fact, since it splits into a product of linear factors over $\F_{q^3}$,
it can be irreducible over $\F_2$ only if its degree $q+1$ divides
$3s$, which occurs only for $s=1,3$.

Of course, many of the elements of $\mathcal{N}_2$ given by Proposition~\ref{prop:2^2s+}
are not really ``new''.
In particular, it is immediate that $q^2+q+1$ is a multiple of $7=2^3-1$
unless $s$ is a multiple of three.
It is also easy to see that
$q^{6}+q^{3}+1$ is a multiple of $73$ unless $s$ is a multiple of three.
In fact, the really ``new'' elements of $\mathcal{N}_2$
produced by Proposition~\ref{prop:2^2s+} are those for which $s$ is a power of three.
We postpone a precise statement and a proof of this fact, in greater generality,
to Proposition~\ref{prop:nearly-coprime}.

The crucial computation in the proof of Proposition~\ref{prop:2^2s+}
generalizes to show the following result, which may be of independent interest.

\begin{prop}\label{prop:matrix}
Let $p$ be a prime, $q=p^s$,
$M=\left[
\begin{array}{cc}
a&b\\
c&d
\end{array}
\right]
\in\mathrm{GL}_2(q)$,
and
$f(x)=cx^{q+1}+dx^q-ax-b\in\F_q[x]$.
Then the period of $f(x)$ divides $q^u-1$, where $u$ is the order of
the image of $M$ in the group $\mathrm{PGL}_2(q)$.
(Equivalently, $f(x)$ splits into a product of linear factors over
$\F_{q^u}$.)

Suppose, in addition, that in the action of $\langle M\rangle$
by right-multiplication on the two dimensional row space over $\F_q$,
the vectors $(1,0)$ and $(0,1)$ belong to the same orbit.
Then the period of $f(x)$ divides $(q^v-1)/(q-1)$,
where $v$ is the order of $M$ in $\mathrm{GL}_2(q)$.
\end{prop}

\begin{proof}
We denote equality of images in the quotient ring
$\F_q[x]/(cx^{q+1}+dx^q-ax-b)$
by a congruence sign.
Since $ad\not=bc$, the binomial $cx+d$ is invertible
in the quotient ring, where we have
$x^q\equiv (ax+b)/(cx+d)$.
By taking the $q$-th powers of both sides we obtain that
$x^{q^2}\equiv (ax^q+b)/(cx^q+d)$,
because $a,b,c,d$ are invariant under
the map $\alpha\mapsto\alpha^q$.
This calls for a further reduction of $x^q$.
More generally, the powers
$x^{q^i}$ (in the quotient ring) can be computed by iterating the
substitution
$x^q\mapsto (ax+b)/(cx+d)$.
Induction shows that
\begin{equation}\label{eq:q-powers}
x^{q^i}
\equiv\frac{ex+f}{gx+h},
\quad\text{where}\quad
\left[
\begin{array}{cc}
e&f\\
g&h
\end{array}
\right]
=
\left[
\begin{array}{cc}
a&b\\
c&d
\end{array}
\right]^i.
\end{equation}
This also follows from the
well-known faithful representation
of $\mathrm{PGL}_2(q)$
as the group of rational expressions of the form
$(ex+f)/(gx+h)\not=1$, with coefficients in $\F_q$,
under substitution.
We conclude that $x^{q^u}\equiv x$, and hence the period of the polynomial
$f(x)$ divides $q^u-1$.

To be able to strengthen our conclusion we impose the additional condition
that the vectors $(1,0)$ and $(0,1)$ belong to the same orbit
in the action of $\langle M\rangle$
by right-multiplication on the two dimensional row space over $\F_q$.
This means that some power of $M$ has $(1,0)$ as its second row.
Let $v$ (a multiple of $u$) be the order of $M$ in $\mathrm{GL}_2(q)$.
Then the numerators of the fractions which can replace the various
factors $x^{q^i}$ in the expression
$x^{(q^v-1)/(q-1)}=
x\cdot x^q\cdot x^{q^2}\cdots x^{q^{v-1}}$
according to formula~\eqref{eq:q-powers},
are the same as as the denominators, just in a different order.
Consequently, they cancel out and we obtain that
$x^{(q^v-1)/(q-1)}\equiv 1$.
Therefore, under the present additional condition, the period of the
polynomial
$cx^{q+1}+dx^q-ax-b$ divides $(q^v-1)/(q-1)$.
\end{proof}

Our proof of Proposition~\ref{prop:2^2s+} is the special case of Proposition~\ref{prop:matrix}
where $q=2^s$ and
$M$ is the matrix
$\left[
\begin{smallmatrix}
1&1\\1&0
\end{smallmatrix}\right]$
of order three.
Polynomials of the form $x^{q+1}-ax-b$ occur in several areas of mathematics and have been extensively studied,
see~\cite{Bluher} and the references therein.
They are a special case of the {\em projective} polynomials of~\cite{Abhyankar},
and their connection with projective linear groups is much deeper than the superficial aspect
employed in the proof of Proposition~\ref{prop:matrix}.
Unfortunately, it does not seem possible to use Proposition~\ref{prop:matrix} to prove that
$(q^t-1)/(q-1)$ belongs to $\mathcal{N}_2$, where $q=2^s$, for
values of $t$ higher than three,
although this conclusion is true, and follows from the more general results of the next section.
See also Proposition~\ref{prop:2^3s+} for a direct proof in case $t=4$.

\section{A Frobenius group}\label{sec:Frobenius}

The main result of this section and of the paper, Theorem~\ref{thm:U-big}, produces many non-trivial elements of $\mathcal{N}_2$.
It shows that any subgroup of $\F_{2^k}^{\ast}$ of order
{\it large enough} with respect to its index contains
an affine $\F_2$-line with direction $1$.
Our proof of Theorem~\ref{thm:U-big} was inspired by an argument
in Chapter~VI of Feit and Thompson's ``Odd Order Paper''~\cite[Lemmas~38.9 and~38.10]{FT},
as simplified in \cite[Lemme~2]{Pet}.
I am grateful to I.~M.~Isaacs for pointing out to me that a generalization of this character-theoretic argument
can be found in~\cite[Section~26]{Feit}.

\begin{theorem}\label{thm:U-big}
Let $q=2^k$ and let $U_n$ be a subgroup of $\F_{q}^{\ast}$ of order $n$, with
$n^4>(q-n)^3$.
Then there exists $\alpha\in U_n$ such that $\alpha+1\in U_n$.
Consequently, $n\in\mathcal{N}_2$.
\end{theorem}

\begin{proof}
Consider the semidirect product $G$ of $\F_{q}$
with an isomorphic copy $\bar U$ of its subgroup $U_n$, acting on
$\F_{q}$ by multiplication.
Let $e$ be the number of elements $\alpha$ of $U_n$ such that $\alpha+1\in U_n$.
If $g$ is any nonzero element of $\F_{q}$,
for example $g=1$, then $e$ is also the number of ordered pairs
$(u,v)$ of elements of $U_n$ such that $gu+gv=g$.
Since the coset $gU_n$ of $U_n$ in $\F_{q}^{\ast}$ coincides with the
conjugacy class  $\mathcal{K}$ of $g$ in $G$,
the number $e$ is the so-called structure
constant of $G$ (strictly speaking, of the center
of the complex group algebra of $G$) with respect to the classes
$\mathcal{K}$, $\mathcal{K}$, $\mathcal{K}$.
Thus $e$ can be computed
in terms of the (complex) characters of $G$,
as in~\cite[Problem (3.9), p.~44]{Isa}.

Since $G$ is a Frobenius group with kernel $\F_{q}$,
it  has $n$ linear characters,
whose kernels contain the derived subgroup $G'=\F_{q}$, and
further irreducible characters $\chi_i$,
for $i=1,\ldots,(q-1)/n$, each of degree $n$.
Noting that $\chi_i(g)$ is an integer and, in particular, is real,
it follows that
\[
e=\frac{\mathcal{|K|}^2}{|G|}
\sum_{\chi\in\Irr(G)}
\frac{\chi(g)^2 \overline{\chi}(g)}{\chi(1)}
=
\frac{n}{q}\left(
 n+
 \frac{1}{n}
 \sum_i \chi_i(g)^3
\right).
\]
Therefore, we have $q e=n^2+\sum_i \chi_i(g)^3$.
If we can show that the absolute value of $\sum_i \chi_i(g)^3$ is less than $n^2$, we conclude that $e>0$.

According to the second orthogonality relation for characters we have
$n+\sum_i \chi_i(g)^2=|\mathbf{C}_G(g)|=q$.
In particular,
we have $|\chi_i(g)|\le (q-n)^{1/2}$ for all $i$, and hence
\[
 \big|
 \sum_i \chi_i(g)^3
 \big|
\le
 \max_i |\chi_i(g)|
 \cdot
 \sum_i \chi_i(g)^2
\le
 (q-n)^{3/2}.
\]
Therefore, under our hypothesis that $n^4>(q-n)^3$ we have
$qe\ge n^2-\big|\sum_i \chi_i(g)^3\big|>0$, which is the desired conclusion.
\end{proof}

A slightly stronger hypothesis on $n$ than that of Theorem~\ref{thm:U-big}, but perhaps easier to remember,
is that the multiplicative group of a field $\F_q$ which contains a subgroup $U_n$ of order $n$ has order
at least the fourth power of the index of $U_n$.

\begin{cor}\label{cor:N_2-examples}
If $q=2^s$ and $t\ge 4$, then $(q^t-1)/(q-1)=q^{(t-1)}+\cdots+q+1$
belongs to $\mathcal{N}_2$.
\end{cor}

\begin{proof}
Set $n=(q^t-1)/(q-1)$.
Because $t\ge 4$ we have $(q-1)^4<q^t-1$, and hence
\[
n^4=(q^t-1)^4/(q-1)^4>(q^t-1)^3>(q^t-n)^3.
\]
Since $n$ divides $q^t-1$ the conclusion follows from Theorem~\ref{thm:U-big}.
\end{proof}

Proposition~\ref{prop:2^2s+} and Corollary~\ref{cor:N_2-examples} together say that
$(q^t-1)/(q-1)$ belongs to $\mathcal{N}_2$ for all $t\ge 3$, where $q=2^s$ as above.
However, the case $t=3$, which we have proved directly in Proposition~\ref{prop:2^2s+},
does not follow from Theorem~\ref{thm:U-big}.
(See Proposition~\ref{prop:2^3s+} for a direct proof in case $t=4$.)
Furthermore, their joint statement does not extend to $t=2$, because
$(q^2-1)/(q-1)=q+1\in\mathcal{N}_2$
if and only if $s$ is odd.
In fact, modulo $x^2-x-c$ we have
$x^{q}\equiv x+c+c^2+c^4+\cdots+c^{q/2}$,
and hence
$x^{q+1}\equiv x^2+x(c+c^2+c^4+\cdots+c^{q/2})
\equiv x(1+c+c^2+\cdots+c^{q/2})+c$.
Since this equals $1$ if and only if $c=1$ and $s$ is odd, our claim follows according to condition~(3)
of Theorem~\ref{thm:n-equiv}.
However, when $s$ is odd the number $q+1$ is a multiple of $3$,
and hence is a trivial element of $\mathcal{N}_2$.

Corollary~\ref{cor:N_2-examples} produces the next smallest non-trivial element of $\mathcal{N}_2$
after $73$, namely, $85=(2^8-1)/(2^2-1)$.
Most of the elements of $\mathcal{N}_2$ given by
Proposition~\ref{prop:2^2s+} and Corollary~\ref{cor:N_2-examples}
are proper multiples of other numbers of the same form.
In Proposition~\ref{prop:nearly-coprime}
we determine those which are not.
Of course, the numbers $(q^t-1)/(q-1)$ with $t\ge 3$ and their multiples do not exhaust $\mathcal{N}_2$, and
we give a few numerical examples in Remark~\ref{rem:examples}.

We will need a few elementary facts about integers of the form $p^a-1$, where $p$ is a prime number.
The simplest is that $p^a-1$ divides $p^b-1$ if and only if $a$ divides $b$.
In fact, this is true if $p$ is any integer different from $0$, $\pm 1$ or $-2$.
However, the case where $p$ is a prime admits a more elegant proof
(leaving aside the trivial case where $a=0$, and hence $b=0$)
by viewing $p^a-1$ as the order
of the multiplicative group of the field $\F_{p^a}$, and noting that $\F_{p^a}$ is a subfield of $\F_{p^b}$
if and only if $a$ divides $b$.
It follows that $(p^a-1,\, p^b-1)=p^{(a,b)}-1$ for any positive integers $a,b$
(where $(a,b)$ denotes the greatest common divisor of $a$ and $b$).
Furthermore, $(p^{ab}-1)/(p^b-1)$ divides $(p^{abc}-1)/(p^{bc}-1)$ if $(a,c)=1$.
In fact, since both $p^{ab}-1$ and $p^{bc}-1$ divide $p^{abc}-1$, and
$(p^{ab}-1,\, p^{bc}-1)=p^b-1$, we have that
$(p^{ab}-1)(p^{bc}-1)$ divides $(p^{abc}-1)(p^b-1)$,
whence the conclusion.
We record the next fact as a lemma.

\begin{lemma}\label{lemma:coprime}
Let $p$ and $r$ be primes (not necessarily distinct).
Then the integers $(p^{r^{a+1}}-1)/(p^{r^a}-1)$
are pairwise coprime, for $a\ge 0.$
\end{lemma}

\begin{proof}
Let $q$ be a prime divisor of
$(p^{r^{a+1}}-1)/(p^{r^a}-1)$.
Then the image of $p^{r^a}$ in $\F_q$ is a root of the polynomial
$(x^r-1)/(x-1)\in\F_q[x]$.
Hence the image of $p$ in $\F_q$ is a nonzero element of multiplicative order (exactly) $r^{a+1}$.
In particular, $q$ determines $a$ uniquely, and the conclusion follows.
\end{proof}

\begin{prop}\label{prop:nearly-coprime}
Every integer of the form
$(q^t-1)/(q-1)$ with $q=2^s$ and $t\ge 3$
is a multiple of at least one element of the set
\[
\mathcal{B}
=
\left\{\left.
\frac{2^{2^{a+2}}-1}{2^{2^a}-1},\,
\frac{2^{r^{b+1}}-1}{2^{r^b}-1}\,
\right\vert\,
r\in\mathbb{P}\setminus\{2\},\,
a\ge 0,\, b\ge 0
\right\},
\]
where $\mathbb{P}$ denotes the set of prime numbers.
The elements of $\mathcal{B}$ are pairwise coprime, with the only exception
of pairs of elements
$(2^{2^{a+2}}-1)/(2^{2^a}-1)$
for consecutive values of $a$, where one has
$\displaystyle\left(
\frac{2^{2^{a+2}}-1}{2^{2^a}-1},\,
\frac{2^{2^{a+3}}-1}{2^{2^{a+1}}-1}
\right)=
\frac{2^{2^{a+2}}-1}{2^{2^{a+1}}-1}
$,
which does not belong to $\mathcal{N}_2$.
In particular, no element of $\mathcal{B}$ is a proper multiple of any number of the form
$(q^t-1)/(q-1)$ with $q=2^s$ and $t\ge 3$.
\end{prop}

\begin{proof}
Consider
$(q^t-1)/(q-1)$,
for some $t\ge 3$.

If $t$ is divisible by an odd prime $r$,
then
$(q^t-1)/(q-1)$
is a multiple of
$(q^r-1)/(q-1)$.
Write $s=r^b\cdot c$ with $c$ prime to $r$.
Then $(q^r-1)/(q-1)$
is a multiple of
$(2^{r^{b+1}}-1)/(2^{r^b}-1)$.

If $t$ is not divisible by an odd prime, then it is a power of two greater than two, and hence
$(q^t-1)/(q-1)$
is a multiple of
$(q^4-1)/(q-1)$.
Write $s=2^a\cdot c$ with $c$ odd.
Then $(q^4-1)/(q-1)$
is a multiple of
$(2^{2^{a+2}}-1)/(2^{2^a}-1)$.

The coprimality statement follows at once from Lemma~\ref{lemma:coprime}.
\end{proof}

Note that the expression $(q^t-1)/(q-1)$ with $q=2^s$ and $t\ge 3$
includes, by taking $s=1$, all integers of the form $2^k-1$ with the exception of $3$.
One may like to include this case to embrace all numbers in $\mathcal{N}_2$ for
which we have found explicit parametric expressions, as follows.
All numbers in $\mathcal{N}_2$ which have some divisor of the form $2^k-1$ or $(q^t-1)/(q-1)$ with $q=2^s$ and $t\ge 3$,
have also a divisor in the subset $\mathcal{B}\cup\{3\}$ of $\mathcal{N}_2$.
In this context we may add that all numbers in $\mathcal{B}$ are prime to $3$,
except for $(2^{2^2}-1)/(2^{2^0}-1)=15$.
In particular, no element of $(\mathcal{B}\setminus\{15\})\cup\{3\}$
is a proper multiple of any other number of the form $2^k-1$ or $(q^t-1)/(q-1)$ with $q=2^s$ and $t\ge 3$.

\begin{cor}\label{cor:n-divisor}
For any positive integer $k$ which is divisible by $8$ or by the square of an odd prime, there is a proper divisor of
$2^k-1$ in $\mathcal{N}_2$ which is not a multiple of any number of the form $2^h-1$.
\end{cor}

\begin{proof}
If $k$ is a multiple of $8$ then $2^k-1$ is a multiple of $85=(2^8-1)/(2^2-1)$.
If $k$ is a multiple of $r^2$ for an odd prime $r$, then $2^k-1$ is a multiple of
$(2^{r^2}-1)/(2^{r}-1)$.
Since the latter belongs to $\mathcal{B}$, it is not a multiple of any number of the form $2^h-1$,
according to Proposition~\ref{prop:nearly-coprime} and the comments which follow it.
\end{proof}

\begin{rem}\label{rem:examples}
Computer calculations based on Lemma~\ref{lemma:remainder}
have shown that the elements of $\mathcal{N}_2$ which are less than $50000$,
which are not proper multiples of other elements of $\mathcal{N}_2$, and are not of the form $2^k-1$, are
$73$, $85$, $3133$,
$4369$,
$11275$ and $49981$.
The first, second and fourth number in this list are predicted by
Proposition~\ref{prop:2^2s+} and Corollary~\ref{cor:N_2-examples}, being
$(2^9-1)/(2^3-1)$,
$(2^8-1)/(2^2-1)$ and
$(2^{16}-1)/(2^4-1)$.
The remaining numbers can be expressed as
$3133=\frac{(2^{24}-1)(2^2-1)}{(2^8-1)(2^6-1)}
=|\F_{2^{24}}^{\ast}:
\langle\F_{2^8}^{\ast},\F_{2^6}^{\ast}\rangle|$,
$11275=\frac{(2^{20}-1)}{(2^5-1)(2^2-1)}
=|\F_{2^{20}}^{\ast}:
\langle\F_{2^5}^{\ast},\F_{2^2}^{\ast}\rangle|$,
and
$49981=\frac{(2^{30}-1)(2^2-1)}{(2^{10}-1)(2^6-1)}
=|\F_{2^{30}}^{\ast}:
\langle\F_{2^{10}}^{\ast},\F_{2^6}^{\ast}\rangle|$.
Since
$3133<(2^{24})^{0.484}$,
$11275<(2^{20})^{0.674}$ and
$49981<(2^{30})^{0.521}$,
these last three numbers are quite far from the range of elements of $\mathcal{N}_2$
produced by Theorem~\ref{thm:U-big}, which are, roughly, the divisors of $2^k-1$ greater than $(2^k)^{0.75}$.
\end{rem}

\section{Counting lines}\label{sec:counting}

We have observed after the proof of Theorem~\ref{thm:n-equiv}
that the positive integers $n$ prime to $p$ which belong to $\mathcal{N}_p$
are those for which $U_n$ contains an affine $\F_p$-line with direction $1$.
More generally, one may ask for the number
of $\F_p$-lines with direction $1$ contained in the subgroup $U_n$
of $\overline{\F}_p^{\ast}$ of order $n$, for specific values of $n$.
The significance of this more general question for Lie algebras, if any, is not clear.
However, it is not unreasonable to expect that posing a more general question
may help to gain a better understanding of the set $\mathcal{N}_p$.
Furthermore, some of the results of this section do have a constructive value for Lie algebras,
which we discuss in Remark~\ref{rem:constructive}.

We introduce some further notation.
For $n$ prime to $p$ we denote by $E_p(n)$ the set of elements
$\alpha\in U_n$ such that $\alpha+\lambda\in U_n$
for all $\lambda\in\F_p$,
and by $e_p(n)$ the cardinality of $E_p(n)$.
The greatest common divisor of the $p$ polynomials
$(x+\lambda)^n-1$ for $\lambda=0,\ldots,p-1$
has exactly the elements of $E_p(n)$ as roots, each with multiplicity one.
In particular, its degree equals $e_p(n)$.
Note that the number of affine $\F_p$-lines with direction $1$ contained in $U_n$ equals $e_p(n)/p$.
For example, since $E_p(p^k-1)=\F_{p^k}\setminus\F_p$ we have $e_p(p^k-1)=p^k-p$.
In the sequel we collect some less trivial cases where we can compute $e_p(n)$
by determining the greatest common divisor of the polynomials
$(x+\lambda)^n-1$ for $\lambda=0,\ldots,p-1$,
thus giving a quite explicit description of the set $E_p(n)$.
We will need the following result, which describes an alternative method for finding the greatest common divisor of the polynomials
$(x+\lambda)^n-1$ for $\lambda=0,\ldots,p-1$ and, in particular, its degree $e_p(n)$.

\begin{lemma}\label{lemma:remainder}
Let $n$ be prime to $p$.
Write the remainder of the division of $x^n-1$ by $x^p-x-c$
(with respect to the indeterminate $x$, and in characteristic $p$)
in the form $r_{p-1}(c)x^{p-1}+\cdots+r_1(c)x+r_0(c)$,
where the coefficients $r_{p-1}(c),\ldots,r_0(c)$ are
polynomials in the indeterminate $c$.
Let $g(c)=(r_{p-1}(c),\ldots,r_0(c))$
be the greatest common divisor of the coefficients.
Then $g(x^p-x)$ equals the greatest common divisor of the polynomials
$(x+\lambda)^n-1$ for $\lambda=0,\ldots,p-1$.
In particular, $e_p(n)$ equals $p$ times the degree of the polynomial $g(c)$.
\end{lemma}

\begin{proof}
Since $\alpha,\alpha+1,\ldots,\alpha+p-1$ are the roots of the
polynomial $x^p-x-(\alpha^p-\alpha)$,
an element $\alpha\in\bar\F_p$ belongs to $E_p(n)$ if and only if
$x^p-x-\bar c$ divides $x^n-1$, where $\bar c=\alpha^p-\alpha$.
This occurs if and only if
$r_{p-1}(\bar c)=\cdots=r_0(\bar c)=0$, or, in turn, if and only if $\bar c$ is a root of $g(c)$.
Thus $E_p(n)$, which is the set of roots of the greatest common divisor of the polynomials
$(x+\lambda)^n-1$ for $\lambda=0,\ldots,p-1$,
is also the set of roots of $g(x^p-x)$.
To complete the proof it remains to show that the polynomial $g(c)$ has no multiple roots.

Write
\[
x^n-1=r_{p-1}(c)x^{p-1}+\cdots+r_1(c)x+r_0(c)+
Q(x,c)(x^p-x-c).
\]
Differentiating with respect to $c$ we obtain
\[
0=r_{p-1}^{\prime}(c)x^{p-1}+\cdots+r_1^{\prime}(c)x+r_0^{\prime}(c)+
Q_c(x,c)(x^p-x-c)-Q(x,c),
\]
where $Q_c$ denotes the partial derivative of $Q$ with respect to $c$.
Suppose for a contradiction that
$g(c)$ has a multiple root $\bar c$.
Then all of $r_{p-1}(c),\ldots,r_0(c)$ have $\bar c$ as a root
with multiplicity greater than one, and so $\bar c$ will be a root
of their derivatives, too.
Substituting $\bar c$ for $c$ in our two equalities we obtain
$x^n-1=Q(x,\bar c)(x^p-x-\bar c)$, and
$0=Q_c(x,\bar c)(x^p-x-\bar c)-Q(x,\bar c)$, whence
$x^n-1=Q_c(x,\bar c)(x^p-x-\bar c)^2$.
This implies that $x^n-1$ has multiple roots, contradicting our
assumption that $n$ is prime to $p$.
\end{proof}

\begin{rem}\label{rem:remainder}
Lemma~\ref{lemma:remainder} implicitly gives a convenient algorithm for computing $e_p(n)$ for a given integer $n$.
In fact, we have used that algorithm to compute tables of elements in $\mathcal{N}_p$
for various small primes $p$, as we have partially reported in Remark~\ref{rem:examples} for $p=2$.
A more obvious algorithm, based on condition~(2) of Theorem~\ref{thm:n-equiv}
rather than on condition~(3),
is computing $e_p(n)$ as the degree of the greatest common divisor of the $p$ polynomials
$(x+\lambda)^n-1$ for $\lambda=0,\ldots,p-1$.
The algorithm given by Lemma~\ref{lemma:remainder}, despite being slightly inferior
for the computation of $e_p(n)$ for a specific value of $n$, performs better than the other
in determining $e_p(n)$ for a long sequence of consecutive values of $n$.
This is essentially because computing the remainder of the division of $x^n-1$ by $x^p-x-c$
is very fast when the remainder of the division of $x^{n-1}-1$ by $x^p-x-c$
is already available, while the other algorithm does not allow such a reduction.
\end{rem}

We now apply Lemma~\ref{lemma:remainder} to refine the statement that $q+1\in\mathcal{N}_2$, where $q=2^s$,
if and only if $s$ is odd, which we have shown in the paragraph following Corollary~\ref{cor:N_2-examples}.
In fact, the same calculation done there
shows that
$g(c)=c-1$ or $1$ (and hence $(x^{q+1}-1,(x+1)^{q+1}-1)=g(x^2+x)$ equals $x^2+x+1$ or $1$)
according as $s$ is odd or even.
Consequently, we have
$e_2(q+1)=1-(-1)^s$.
Now we use Lemma~\ref{lemma:remainder} to compute $e_2(n)$
for $n=(q^3-1)/(q-1)=q^2+q+1$.

\begin{prop}\label{prop:e-values}
Let $q=2^s$.
Then the greatest common divisor
$(x^{q^2+q+1}-1,(x+1)^{q^2+q+1}-1)$
equals
$(x^{q+1}+x+1)(x^{q+1}+x^{q}+1)$ for $s$ odd and
$(x^{q+1}+x+1)(x^{q+1}+x^{q}+1)/(x^2+x+1)$ for $s$ even.
In particular, we have
$e_2(q^2+q+1)=2q+1-(-1)^s$.
\end{prop}

\begin{proof}
Let $n=q^2+q+1$.
All congruences will tacitly be modulo $x^2+x+c$.
An easy induction, already employed elsewhere, shows that
$x^{q}\equiv x+\gamma$, where
$\gamma=\sum_{i=0}^{s-1} c^{2^i}$, and so
$x^{q^2}\equiv (x+\gamma)^{q}=x+\gamma+\gamma^{q}$.
It follows that
\begin{eqnarray*}
x^{q^2+q+1}
&\equiv&
(x+\gamma+\gamma^{q})(x+\gamma)x
\equiv
(x+\gamma+\gamma^{q})(x(1+\gamma)+c) \\
&\equiv&
x^2(1+\gamma)+x((\gamma+\gamma^{q})(1+\gamma)+c)
 +(\gamma+\gamma^{q})c \\
&\equiv&
x((1+\gamma+\gamma^{q})(1+\gamma)+c)+(1+\gamma^{q})c.
\end{eqnarray*}
Hence, we have $x^{n}-1\equiv r_1(c)x+r_0(c)$, where
$r_1(c)=(1+\gamma+\gamma^{q})(1+\gamma)+c$ and
$r_0(c)=(1+\gamma^{q})c+1$.
In order to apply Lemma~\ref{lemma:remainder} we need to compute $g(c)=(r_1(c),\, r_0(c))$.

Since $\gamma^2+\gamma=c^{q}+c$ we have
\[
r_1(c)\cdot c+r_0(c)\cdot (\gamma+1)=
(\gamma(1+\gamma)+c)c+(1+\gamma)=
c^{q+1}+\gamma+1.
\]
Therefore, $g(c)$
divides $c^{q+1}+\gamma+1$.
Actually, the derivative criterion shows that $1$ is a multiple root of
$c^{q+1}+\gamma+1$ when $s$ is even.
Since $g(c)$ has no multiple roots,
it must divide $(c^{q+1}+\gamma+1)/(c+1)$ in that case.
(Incidentally, we know that $1$ must be a root of $g(c)$ when $s$ is even
because $x^2+x+1$, which divides $x^3-1$, divides also $x^{q^2+q+1}-1$ in that case.)
According to Lemma~\ref{lemma:remainder}, it follows that
$e_2(q^2+q+1)=2\cdot\deg(g(c))\le 2q+1-(-1)^s$.

Now we prove the reverse inequality.
In the proof of Proposition~\ref{prop:2^2s+} we have shown that any root $\alpha$ of
the polynomial $x^{q+1}+x+1$ satisfies
$\alpha^{q^2+q+1}=1$ and
$(\alpha+1)^{q^2+q+1}=1$.
Note also that if any $\alpha$ satisfies both these conditions,
then $(\alpha+1)/\alpha=1+1/\alpha$ also does.
Therefore, the roots of the reciprocal polynomial
$x^{q+1}+x^{q}+1$ of $x^{q+1}+x+1$ also satisfy
$\{\alpha,\, \alpha+1\}\subseteq U_{q^2+q+1}$, and so do
all roots of the product $(x^{q+1}+x+1)(x^{q+1}+x^{q}+1)$.
Since each of the two factors has no multiple roots, and their greatest common divisor
is or $1$ or $x^2+x+1$ according as $s$ is odd or even, we conclude that
$|E_2(q^2+q+1)|\ge 2q+1-(-1)^s$.

Taking into account the first part of the proof, equality holds here.
It follows that
$g(x^2+x)$ equals $(x^{q+1}+x+1)(x^{q+1}+x^{q}+1)$ for $s$ odd and
$(x^{q+1}+x+1)(x^{q+1}+x^{q}+1)/(x^2+x+1)$ for $s$ even.
\end{proof}

The fact that the polynomial $g(x^2-x)=(x^n-1,(x+1)^n-1)$ is divisible by such a nice trinomial as $x^{q+1}+x+1$
is the reason why the case $t=3$ of $n=(q^t-1)/(q-1)$
can be dealt with so explicitly in Propositions~\ref{prop:2^2s+} and~\ref{prop:e-values}.
Such explicitness seems not easy to achieve when $t>3$.
In particular, already when $t=4$ one has to work harder just to prove that $g(x^2-x)$ is not constant,
as in the proof of the following result.
As noted in the Introduction, Proposition~\ref{prop:2^3s+} is also a consequence of Corollary~\ref{cor:N_2-examples},
but its direct proof given here has an added value which we point out in Remark~\ref{rem:constructive}.

\begin{prop}\label{prop:2^3s+}
If $q=2^s$, then the number
$(q^4-1)/(q-1)=q^3+q^2+q+1$
belongs to $\mathcal{N}_2$.
\end{prop}

\begin{proof}
Consider the polynomial $x^{q^2}+x+1$ over $\F_2$. Any root $\alpha$ of this polynomial in $\bar\F_2$ is Galois
conjugate to $\alpha+1$, because $\alpha^{q^2}=\alpha+1$.
Therefore, every divisor of $x^{q^2}+x+1$ over
$\F_2$ is invariant under the substitution $x\mapsto x+1$.
Furthermore, the period of $x^{q^2}+x+1$ divides $q^4-1$, because
$\alpha^{q^4}=(\alpha+1)^{q^2}=\alpha^{q^2}+1=\alpha$.
Since
\[
x^{q^3+q^2+q+1}=(x^{q^2+1})^{q+1}
\equiv ((x+1)x)^{q+1}
\not\equiv 1\pmod{x^{q^2}+x+1},
\]
the period of $x^{q^2}+x+1$ does not divide $q^3+q^2+q+1$.
However, the period of the greatest common divisor $h(x)$ of $(x^2+x)^{q+1}-1$ and $x^{q^2}+x+1$
does divide $q^3+q^2+q+1$.
Since we have seen that $h(x)$ must be invariant under the substitution $x\mapsto x+1$,
for any of its roots $\alpha$ the orders of both $\alpha$ and $\alpha+1$ divide
$q^3+q^2+q+1$.
In order to conclude that this number belongs to $\mathcal{N}_2$ it remains to show that $h(x)$ is not a constant polynomial.

Since the polynomials $(x^2+x)^{q+1}-1$ and $x^{q^2}+x+1$
can be obtained by substituting $x^2+x$ for $y$ in
$y^{q+1}-1$ and $y^{q^2/2}+y^{q^2/4}+\cdots+y^2+y+1$,
it is enough to prove that these two polynomials in the indeterminate $y$
have a nonconstant common factor.
The latter polynomial divides $y^{q^2}-y$, and hence $\F_{q^2}$ contains a splitting field for it.
In fact, the roots of that polynomials consist of all elements of $\F_{q^2}$ of absolute trace $1$.
The roots of the former polynomial form the subgroup of $\F_{q^2}^\ast$ of order $q+1$,
and hence span $\F_{q^2}$ over $\F_2$.
(Any multiplicative subgroup spans a subfield, which must coincide with the whole field in this case.)
Consequently, they cannot be all contained in the $\F_2$-hyperplane of $\F_{q^2}$ consisting
of the elements of absolute trace zero.
We conclude that the two polynomials have a common root in $\F_{q^2}$, and hence a common factor over $\F_2$.
\end{proof}

\begin{rem}\label{rem:constructive}
The construction of a non-nilpotent Lie algebra with a nonsingular derivation of order $n$
given in~\cite[Theorem~2.1]{Mat:nonsing-der} requires an element $\alpha\in E_p(n)$, that is,
a common root of the polynomials
$(x+\lambda)^n-1$ for $\lambda=0,\ldots,p-1$.
In order to explicitly construct an admissible element $\alpha$ of $\bar\F_p$ in specific instances,
it is therefore of interest to produce a common divisor of these polynomials over the prime field $\F_p$,
having relatively low degree with respect to $n$.
The proofs of Propositions~\ref{prop:2^2s+} and~\ref{prop:2^3s+} are both based on producing such a common divisor, namely,
$x^{q+1}+x+1$ in the former case
and the polynomial denoted by $h(x)$ in the latter.
In the latter case, we have given $h(x)$ as the greatest common divisor
of two further polynomials, of degrees $2q+2$ and $q^2$,
which are both much smaller than $n$ when $q$ is large.
In this sense, the proofs of Propositions~\ref{prop:2^2s+} and~\ref{prop:2^3s+}
have a constructive character which is missing in the method of Section~\ref{sec:Frobenius}.
We have been unable to explicitly produce such a common divisor (and thus a proof analogous to those just mentioned) in case
$n=(q^t-1)/(q-1)$ with $t>4$.
\end{rem}

We conclude this section by computing $e_p(n)$ for $n=(p^p-1)/(p-1)$,
thus providing a proof for an equivalent statement mentioned at the end of Section~3 of~\cite{Mat:nonsing-der}.

\begin{prop}\label{prop:norm-one}
Let $n=(p^p-1)/(p-1)$.
Then the greatest common divisor of the polynomials
$(x+\lambda)^n-1$ for $\lambda=0,\ldots,p-1$
equals $x^p-x-1$.
Consequently, we have $e_p(n)=p$.
\end{prop}

\begin{proof}
We prepare for an application of Lemma~\ref{lemma:remainder}.
Induction shows that
\[
x^{p^i}\equiv
x+\sum_{j=0}^{i-1}c^{p^j}
\pmod{x^p-x-c}
\]
for all $i\ge 0$.
It follows that
$x^n-1=x^{1+p+p^2+\cdots+p^{p-1}}-1$ is congruent to
\begin{equation*}
h(x,c):=x(x+c)(x+c+c^p)
\cdots
(x+c+c^p+\cdots+c^{p^{p-2}})-1
\end{equation*}
modulo $x^p-x-c$.
Since $h(x)$ has leading term $x^p$, the remainder of the division
of $x^n-1$ by $x^p-x-c$ required by
Lemma~\ref{lemma:remainder} equals $h(x,c)-(x^p-x-c)$.
It follows that $r_0(c)=c-1$, and hence the greatest common divisor $g(c)$ of the coefficients of the remainder divides $c-1$.
However, the remainder vanishes for $c=1$, because $h(x,1)=x^p-x-1$,
and we conclude that $g(c)=c-1$.
\end{proof}

It is worth noting that $n=(p^p-1)/(p-1)$ is one case where $e_p(n)$ assumes its smallest possible positive value.
With an imprecise but perhaps suggestive phrasing, we may say that
``$n=(p^p-1)/(p-1)$ does belong to $\mathcal{N}_p$, but just barely''.

\begin{rem}
The full set of affine $\F_p$-lines contained in the multiplicative subgroup $U_n$ of $\F_{p^p}$ with $n=(p^p-1)/(p-1)$
forms an interesting configuration.
In fact, it follows easily from Proposition~\ref{prop:norm-one}
that $U_n$ contains exactly $n$ affine $\F_p$-lines, one for every possible direction.
More precisely, for every $\beta\in U_n$ the unique affine $\F_p$-line in $U_n$ with direction $\beta$
consists of the roots of the polynomial
$x^p-\beta^{p-1}x-\beta^{p}$.
This also implies that each element of $U_n$ belongs to exactly $p$ of these lines.
\end{rem}

\bibliography{References}

\def\cprime{$'$} \def\cprime{$'$}
\providecommand{\bysame}{\leavevmode\hbox to3em{\hrulefill}\thinspace}
\providecommand{\MR}{\relax\ifhmode\unskip\space\fi MR }
\providecommand{\MRhref}[2]{%
  \href{http://www.ams.org/mathscinet-getitem?mr=#1}{#2}
}
\providecommand{\href}[2]{#2}
\begin{thebibliography}{CMN97}

\bibitem[Abh97]{Abhyankar}
Shreeram~S. Abhyankar, \emph{Projective polynomials}, Proc. Amer. Math. Soc.
  \textbf{125} (1997), no.~6, 1643--1650. \MR{MR1403111 (98a:12001)}

\bibitem[AF55]{AF}
A.~A. Albert and M.~S. Frank, \emph{Simple {L}ie algebras of characteristic
  {$p$}}, Univ. e Politec. Torino. Rend. Sem. Mat. \textbf{14} (1954--55),
  117--139. \MR{MR0079222 (18,52a)}

\bibitem[Bar00]{Bartholdy}
Laurent Bartholdi, \emph{Lamps, factorizations, and finite fields}, Amer. Math.
  Monthly \textbf{107} (2000), no.~5, 429--436. \MR{MR1763395 (2001e:11118)}

\bibitem[BKK95]{BKK}
Georgia Benkart, Alexei~I. Kostrikin, and Michael~I. Kuznetsov,
  \emph{Finite-dimensional simple {L}ie algebras with a nonsingular
  derivation}, J. Algebra \textbf{171} (1995), no.~3, 894--916. \MR{MR1315926
  (96b:17020)}

\bibitem[Blu04]{Bluher}
Antonia~W. Bluher, \emph{On {$x\sp {q+1}+ax+b$}}, Finite Fields Appl.
  \textbf{10} (2004), no.~3, 285--305. \MR{MR2067599 (2005b:12005)}

\bibitem[CM05]{CaMa:Hamiltonian}
A.~Caranti and S.~Mattarei, \emph{Gradings of non-graded {H}amiltonian {L}ie
  algebras}, J. Austral. Math. Soc. Ser. A \textbf{79} (2005), no.~3, 399--440.

\bibitem[CMN97]{CMN}
A.~Caranti, S.~Mattarei, and M.~F. Newman, \emph{Graded {L}ie algebras of
  maximal class}, Trans. Amer. Math. Soc. \textbf{349} (1997), no.~10,
  4021--4051. \MR{MR1443190 (98a:17027)}

\bibitem[CN00]{CN}
A.~Caranti and M.~F. Newman, \emph{Graded {L}ie algebras of maximal class.
  {II}}, J. Algebra \textbf{229} (2000), no.~2, 750--784. \MR{MR1769297
  (2001g:17041)}

\bibitem[Fei67]{Feit}
Walter Feit, \emph{Characters of finite groups}, W. A. Benjamin, Inc., New
  York-Amsterdam, 1967. \MR{MR0219636 (36 \#2715)}

\bibitem[FT63]{FT}
Walter Feit and John~G. Thompson, \emph{Solvability of groups of odd order},
  Pacific J. Math. \textbf{13} (1963), 775--1029. \MR{MR0166261 (29 \#3538)}

\bibitem[Isa94]{Isa}
I.~Martin Isaacs, \emph{Character theory of finite groups}, Dover Publications
  Inc., New York, 1994, Corrected reprint of the 1976 original [Academic Press,
  New York]. \MR{MR1280461}

\bibitem[Jac55]{Jacobson:autom-and-deriv}
N.~Jacobson, \emph{A note on automorphisms and derivations of {L}ie algebras},
  Proc. Amer. Math. Soc. \textbf{6} (1955), 281--283. \MR{MR0068532 (16,897e)}

\bibitem[Jac79]{Jac:Lie_algebras}
Nathan Jacobson, \emph{Lie algebras}, Dover Publications Inc., New York, 1979,
  Republication of the 1962 original. \MR{MR559927 (80k:17001)}

\bibitem[LGM02]{L-GMcKay}
C.~R. Leedham-Green and S.~McKay, \emph{The structure of groups of prime power
  order}, London Mathematical Society Monographs. New Series, vol.~27, Oxford
  University Press, Oxford, 2002, Oxford Science Publications. \MR{MR1918951
  (2003f:20028)}

\bibitem[LGN80]{L-GN}
C.~R. Leedham-Green and M.~F. Newman, \emph{Space groups and groups of
  prime-power order. {I}}, Arch. Math. (Basel) \textbf{35} (1980), no.~3,
  193--202. \MR{MR583590 (81m:20029)}

\bibitem[LN83]{LN}
Rudolf Lidl and Harald Niederreiter, \emph{Finite fields}, Encyclopedia of
  Mathematics and its Applications, vol.~20, Addison-Wesley Publishing Company
  Advanced Book Program, Reading, MA, 1983, With a foreword by P. M. Cohn.
  \MR{MR746963 (86c:11106)}

\bibitem[Mat02]{Mat:nonsing-der}
S.~Mattarei, \emph{The orders of nonsingular derivations of modular {L}ie
  algebras}, Israel J. Math. \textbf{132} (2002), 265--275. \MR{MR1952625
  (2003k:17024)}

\bibitem[Pet84]{Pet}
Thomas Peterfalvi, \emph{Simplification du chapitre {VI} de l'article de {F}eit
  et {T}hompson sur les groupes d'ordre impair}, C. R. Acad. Sci. Paris S\'er.
  I Math. \textbf{299} (1984), no.~12, 531--534. \MR{MR770439 (86d:20020)}

\bibitem[Sha94a]{Sha:max}
Aner Shalev, \emph{Simple {L}ie algebras and {L}ie algebras of maximal class},
  Arch. Math. (Basel) \textbf{63} (1994), no.~4, 297--301. \MR{MR1290602
  (95j:17025)}

\bibitem[Sha94b]{Sha:coclass}
\bysame, \emph{The structure of finite {$p$}-groups: effective proof of the
  coclass conjectures}, Invent. Math. \textbf{115} (1994), no.~2, 315--345.
  \MR{MR1258908 (95j:20022b)}

\bibitem[Sha99]{Sha:nonsing-der}
\bysame, \emph{The orders of nonsingular derivations}, J. Austral. Math. Soc.
  Ser. A \textbf{67} (1999), no.~2, 254--260, Group theory. \MR{MR1717417
  (2000k:17021)}

\bibitem[Str04]{Strade:book}
Helmut Strade, \emph{Simple {L}ie algebras over fields of positive
  characteristic. {I}}, de Gruyter Expositions in Mathematics, vol.~38, Walter
  de Gruyter \& Co., Berlin, 2004, Structure theory. \MR{MR2059133
  (2005c:17025)}

\end{thebibliography}

\end{document}